# On the bifurcation sets of functions definable in o–minimal structures

Ta Lê Loi and Alexandru Zaharia

## 1   Introduction

**1.1** Let $g : X \longrightarrow Y$ be a smooth ( i.e. $C^\infty$ differentiable ) map between two smooth manifolds . In analogy with the case of complex polynomial functions , we say that $y_0 \in Y$ is a *typical value* of $g$ if there exists an open neighbourhood $U$ of $y_0$ in $Y$ , such that the restriction $g : g^{-1}(U) \longrightarrow U$ is a $C^\infty$ trivial fibration . If $y_0 \in Y$ is not a typical value of $g$ , then $y_0$ is called an *atypical value* of $g$ . We denote by $\mathrm{B}_g$ the *bifurcation set* of $g$ , i.e. the set of atypical values of $g$ . In the case of a complex polynomial function $f : \mathbf{C}^n \longrightarrow \mathbf{C}$ it is known that $\mathrm{B}_f$ is a finite set , see for instance [12] , [10] . In [1] it is proved that the bifurcation sets of real polynomial functions are also finite . See also [2] for a more general result .

The aim of this note is to show that the bifurcation set $\mathrm{B}_f$ of a smooth *definable* function ( see 1.2 ) $f : \mathbf{R}^n \longrightarrow \mathbf{R}$ is finite . We proceed as in [8] and we estimate from above the bifurcation set . As a by–product of our method , we present a new proof for the fact that the bifurcation set $\mathrm{B}_f$ of a complex polynomial function $f : \mathbf{C}^n \longrightarrow \mathbf{C}$ is finite . Before stating our results , we give some definitions .

**1.2 Definition .** An *o–minimal structure* on the real field $(\mathbf{R}, +, \cdot)$ is a family $\mathcal{D} = (\mathcal{D}_n)_{n \in \mathbf{N}}$ such that , for each $n \in \mathbf{N}$ :

- (D1) $\mathcal{D}_n$ is a boolean algebra of subsets of $\mathbf{R}^n$ , i.e. $\mathcal{D}_n$ is closed under taking the complement and under finite unions .
- (D2) If $A \in \mathcal{D}_n$ , then $A \times \mathbf{R}$ and $\mathbf{R} \times A \in \mathcal{D}_{n+1}$ .
- (D3) If $A \in \mathcal{D}_{n+1}$ , then $\pi(A) \in \mathcal{D}_n$ , where $\pi : \mathbf{R}^{n+1} \longrightarrow \mathbf{R}^n$ is the projection on the first $n$ coordinates .
- (D4) $\mathcal{D}_n$ contains $\{x \in \mathbf{R}^n \mid P(x) = 0\}$ for every polynomial $P \in \mathbf{R}[X_1, \cdots, X_n]$ .
- (D5) Each set belonging to $\mathcal{D}_1$ is a finite union of intervals and points . ( This property is called *o–minimality* . )

A subset $A \subseteq \mathbf{R}^n$ will be called *definable* ( in the structure $\mathcal{D}$ ) if $A \in \mathcal{D}_n$ . A map $f : A \longrightarrow B$ will be called *definable* if the graph of $f$ is a definable set .

O–minimal structures have many nice properties similar to those of semi–algebraic sets . We mention here only the fact that any definable set has only finitely many connected components , and each of them is path connected and definable . For more details , see for instance [3] and [4] . See also [11] .

**1.3** Let $g : \mathbf{R}^n \longrightarrow \mathbf{R}$ be a smooth function . We denote :

$$M(g) := \{\, x \in \mathbf{R}^n \mid \exists \lambda \in \mathbf{R},\ \operatorname{grad} g(x) = \lambda x \,\} \quad .$$

Note that $u \in M(g) \setminus \{0\}$ if and only if either $u$ is a critical point of $g$ , or $u$ is not a critical point of $g$ and the level set $g^{-1}(g(u))$ is a smooth submanifold of $\mathbf{R}^n$ near $u$ , which is not transversal to the sphere $\{\, x \in \mathbf{R}^n \mid \|x\| = \|u\| \,\}$ at $u$ . For a sequence $\{y^k\} \subseteq M(g)$ we consider the conditions

$$\lim_{k \to \infty} \|y^k\| = \infty \ \text{ and } \ \lim_{k \to \infty} g(y^k) = c \quad . \tag{1}$$

We denote by $\Sigma_g$ the set of critical values of $g$ , and we put

$$S_g := \{\, c \in \mathbf{R} \mid \text{there exists a sequence } \{y^k\} \subseteq M(g) \text{ such that (1) is fulfilled} \,\} \quad .$$

The motivation for considering the set $S_g$ is given by the following result , which is similar to Theorem 1 in [8] .

**1.4 Proposition .** *Let $g : \mathbf{R}^n \longrightarrow \mathbf{R}$ be a $C^\infty$ differentiable function . Then for any open interval $U \subseteq g(\mathbf{R}^n) \setminus (\Sigma_g \cup S_g)$ , the restriction*

$$g : g^{-1}(U) \longrightarrow U$$

*is a $C^\infty$ trivial fibration .*

**1.5** Let $\mathcal{D}$ be a fixed , but arbitrary , o–minimal structure on $(\mathbf{R}, +, \cdot)$ . "Definable" will mean definable in $\mathcal{D}$ .

**Theorem .** *Let $f : \mathbf{R}^n \longrightarrow \mathbf{R}$ be a smooth definable function . Then $\Sigma_f$ and $S_f$ are finite sets .*

As a consequence of Proposition 1.4 and of Theorem 1.5 , we have

**Corollary .** *Let $f : \mathbf{R}^n \longrightarrow \mathbf{R}$ be a definable function of class $C^\infty$ . Then $B_f$ is a finite set .*

**1.6** Let $\mathbf{S}_R := \{\, x \in \mathbf{R}^n \mid \|x\| = R \,\}$ . Then the cardinality of $\mathrm{B}_f$ can be estimated from above by using the following

**Proposition .** *Let $f : \mathbf{R}^n \longrightarrow \mathbf{R}$ be a smooth definable function and let $R \in \mathbf{R}$ be a sufficiently large positive number . Then the cardinality of $\Sigma_f$ is less than or equal to the number of connected components of $\{x \in \mathbf{R}^n \mid \operatorname{grad} f(x) = 0\}$ , and the cardinality of $S_f$ is less than or equal to the number of connected components of the intersection $M(f) \cap \mathbf{S}_R$ .*

**1.7** In the complex case , the gradient of a complex polynomial function $f : \mathbf{C}^n \longrightarrow \mathbf{C}$ is defined to be

$$\operatorname{grad} f(x) := \left( \overline{\frac{\partial f}{\partial x_1}(x)}, \ldots, \overline{\frac{\partial f}{\partial x_n}(x)} \right) \quad,$$

where $x = (x_1, \cdots, x_n) \in \mathbf{C}^n$ , the bar denotes the conjugation . We also define

$$M(f) := \{\, x \in \mathbf{C}^n \mid \exists\, \lambda \in \mathbf{C} \ , \ \operatorname{grad} f(x) = \lambda x \,\}$$

and

$$S_f := \{\, c \in \mathbf{C} \mid \text{there exists a sequence } \{y^k\} \subseteq M(f) \text{ such that (1) is fulfilled}\,\} \ .$$

**Theorem .** *Let* $f : \mathbf{C}^n \longrightarrow \mathbf{C}$ *be a complex polynomial function . Then* $S_f$ *is a finite set . Consequently , the bifurcation set* $B_f \subseteq \Sigma_f \cup S_f$ *is finite .*

The next section contains the proofs . Several examples and related remarks are given in the last section.

## 2  Proofs

**2.1** We keep the notations from Section 1 . We denote by $\langle a, b \rangle$ the Euclidean scalar product of $a, b \in \mathbf{R}^n$ . The following Lemma is a direct consequence of the definition .

**Lemma .** *Let* $g : \mathbf{R}^n \longrightarrow \mathbf{R}$ *be a smooth function and let* $J \subseteq \mathbf{R} \setminus S_g$ *be a compact interval . Then* $g^{-1}(J) \cap M(g)$ *is bounded .*

**2.2 Proof of Proposition 1.4 .** We follow closely the proof of Theorem 1 in [8] . Let $U \subseteq g(\mathbf{R}) \setminus (\Sigma_g \cup S_g)$ be an open interval and let $c \in U$ be fixed . From Lemma 2.1 it follows that for any compact subinterval $J \subseteq U$ with $c \in J$ , there exists a sufficiently large number $R \in (0, \infty)$ such that

$$g^{-1}(J) \cap M(g) \cap \{\, x \in \mathbf{R}^n \mid \|x\| \geq R \,\} = \emptyset \quad.$$

Hence for any $x \in A := g^{-1}(J) \cap \{\, x \in \mathbf{R}^n \mid \|x\| \geq R \,\}$ the vectors $\operatorname{grad} g(x)$ and $x$ are linearly independent . Therefore we can find a smooth vector field $v_1(x)$ , defined on $A$ , such that $\langle v_1(x), \operatorname{grad} g(x) \rangle = 1$ and $\langle v_1(x), x \rangle = 0$ . By integrating this vector field one can obtain a $C^\infty$ trivialization for the restriction $g : A \longrightarrow J$ .

Since $J \cap \Sigma_g = \emptyset$ , by Ehresmann's fibration theorem , the restriction $g : g^{-1}(J) \cap \{\, x \in \mathbf{R}^n \mid \|x\| \leq R \,\} \longrightarrow J$ is a $C^\infty$ trivial fibration . Proposition 1.4 follows by glueing together these two trivializations .

$\diamondsuit$

**2.3 Proof of Theorem 1.5 .** By Cell Decomposition Theorem 4.2 in [4] , there exists a partition $\mathcal{S}$ of $\mathbf{R}^n$ into finitely many connected $C^1-$submanifolds , compatible with $\{\ x \in \mathbf{R}^n \ |\ \operatorname{grad} f(x) = 0\ \}$ . Note that if $\Gamma \in \mathcal{S}$ is such that $(\operatorname{grad} f)|_\Gamma = 0$ , then $f|_\Gamma = \operatorname{const}$ . Hence $\Sigma_f = \{\ f(x)\ |\ x \in \Gamma\ ,\ \Gamma \in \mathcal{S}\ ,\ \operatorname{grad} f(x) = 0\}$ is a finite set .

To show that $S_f$ is a finite set , consider

$$A := \{\ (c, y, \lambda, t) \in \mathbf{R} \times \mathbf{R}^n \times \mathbf{R} \times \mathbf{R}\ |\ c = f(y)\ ,\ \operatorname{grad} f(y) = \lambda y\ ,\ t\,\|y\| > 1\ \}\ .$$

Then $A$ is a definable set . Note that

$$c \in S_f \iff \forall \varepsilon > 0\ \ \forall \delta > 0\ \ \exists (c', y, \lambda, t) \in A\ ,\ |t| < \varepsilon\ ,\ |c' - c| < \delta\ .$$

Using the interpretation of the logical symbols in terms of operations on sets , one can see that $S_f$ is definable . Hence $S_f$ is a finite union of points and intervals .

Suppose to the contrary that $S_f$ contains an interval . Let $\pi : \mathbf{R} \times \mathbf{R}^{n+2} \longrightarrow \mathbf{R}$ be the projection on the first coordinate . Since $\pi(A)$ is a finite union of points and intervals and since $S_f \subseteq \overline{\pi(A)}$ , there exists an interval $(\alpha, \beta)$ contained in $S_f \cap \pi(A)$ . Therefore , for all $c \in (\alpha, \beta)$ and for all $t > 0$ , there exists $y(c, t) \in \mathbf{R}^n$ and $\lambda(c, t) \in \mathbf{R}$ such that

$$f(y(c, t)) = c \tag{2}$$

$$\operatorname{grad} f(y(c, t)) = \lambda(c, t)\, y(c, t) \tag{3}$$

$$\|y(c, t)\| > \frac{1}{t} \tag{4}$$

Since $\Sigma$ is a finite set , we may assume , after eventually shrinking the interval $(\alpha, \beta)$ , that $\lambda(c, t) \neq 0$ for all $c \in (\alpha, \beta)$ and for all $t > 0$ . Moreover , by Definable Choice Theorem 4.5 and Cell Decomposition Theorem 4.2 in [4] , there exist a subinterval $(\alpha', \beta') \subseteq (\alpha, \beta)$ and $\varepsilon > 0$ , such that the function $(\alpha', \beta') \times (0, \varepsilon) \ni (c, t) \longmapsto y(c, t)$ is of class $C^1$ . Differentiating with respect to $t$ the relation (2) , we obtain

$$\left\langle \operatorname{grad} f(y(c, t)), \frac{\partial y}{\partial t}(c, t) \right\rangle = 0\ ,\quad \text{for all}\ \ (c, t) \in (\alpha', \beta') \times (0, \varepsilon)\ .$$

Combining this relation with (3) , we get

$$\frac{\partial \|y\|^2}{\partial t}(c, t) = 2 \left\langle y(c, t), \frac{\partial y}{\partial t}(c, t) \right\rangle = 0\ ,\quad \text{for all}\ \ (c, t) \in (\alpha', \beta') \times (0, \varepsilon)\ .$$

Hence , $\|y(c, t)\|$ is independent of $t$ , which contradicts (4) .

$\diamondsuit$

**2.4 Proof of Proposition 1.6 .** The estimation for the cardinality of $\Sigma_f$ is obvious . For the second estimation , we note that the function $\rho : M(f) \longrightarrow \mathbf{R}$ ,

$\rho(x) := \|x\|^2$ is definable . Therefore , by van den Dries' version of Hardt's Theorem on trivialization for definable functions , see [3] , there exists a finite subset $B \subseteq \mathbf{R}$ such that above any connected component of $\mathbf{R} \setminus B$ , the function $\rho$ is a topologically trivial fibration . Hence , there exists $R_0 \in \mathbf{R}$ such that if $R \geq R_0$ , then $M(f) \cap \{x| \ \|x\| \geq R\}$ and $M(f) \cap \mathbf{S}_R$ have the same number of connected components , number which does not depend on $R$ .

Since $S_f$ is finite and by the definition of $S_f$ , it is easy to see that the cardinal of $S_f$ is less than or equal to the number of connected components of $M(f) \cap \{x| \ \|x\| \geq R\}$ , $R \geq R_0$ . The proposition follows .

$\diamondsuit$

**2.5** The proof of Theorem 1.7 will use the following lemmas

**Lemma .** *Let $U$ be an open definable subset of $\mathbf{R}^k$ and let $F : U \times (0, \varepsilon) \longrightarrow \mathbf{R}^m$ be a $C^1$ definable map . Suppose that there exists a constant $K > 0$ such that $\|F(s,t)\| \leq K$ for all $(s,t) \in U \times (0, \varepsilon)$ . Then there exist a definable set $V$ , closed in $U$ and with $\dim V < \dim U$ , and continuous definable functions $\kappa, \tau : U \setminus V \longrightarrow (0, \infty)$ , such that for all $s \in U \setminus V$ and $t \in (0, \tau(s))$ :*

$$\left\| \frac{\partial F}{\partial s}(s,t) \right\| \leq \kappa(s) \quad .$$

For the proof of this lemma , see [5] .

**2.6 Lemma .** *Let $F : (0, \varepsilon) \longrightarrow \mathbf{R}$ be a $C^1$ definable function with $\lim_{t \to 0} F(t) = 0$ . Then $\lim_{t \to 0} tF'(t) = 0$ .*

**Proof** By Monotonicity Theorem 6.1 in [4] , $F$ is either constant , or strictly monotone near 0 . So , it is sufficient to consider the case when $F$ and $F'$ are strictly monotone on $(0, \varepsilon)$ and $F > 0$ . Then $F' > 0$ . By Mean Value Theorem , we have $F(t) = F'(\zeta(t))t$ for $\zeta(t) \in (0,t)$ . It is easy to see that $\zeta : (0, \varepsilon) \longrightarrow (0, \varepsilon)$ is definable and $\lim_{t \to 0} \zeta(t) = 0$ . Therefore ,

$$0 \leq \lim_{\zeta \to 0} \zeta F'(\zeta) = \lim_{t \to 0} \zeta(t) F'(\zeta(t)) \leq \lim_{t \to 0} tF'(\zeta(t)) = 0 \quad .$$

$\diamondsuit$

**2.7 Proof of Theorem 1.7 .** In this proof , we will denote by $\langle a, b \rangle$ the Hermitian product of $a, b \in \mathbf{C}^n$ .

Since the class of all semi–algebraic sets is an o–minimal structure on $(\mathbf{R}, +, \cdot)$ , it follows , as in the proof of Theorem 1.5 , that $S_f$ is definable , i.e. a semi–algebraic subset of $\mathbf{R}^2 \cong \mathbf{C}$ . We will show that $\dim S_f \leq 0$ .

If $\dim S_f \geq 1$, then there exist an interval $J \subseteq \mathbf{R}$ and a smooth curve $\ell : J \longrightarrow S_f$ whose derivative satisfies $\ell'(s) \neq 0$ , for all $s \in J$ . Therefore , for all $s \in J$ and $t > 0$ there exist $p(s,t) \in \mathbf{C}^n$ and $\lambda(s,t) \in \mathbf{C}$ such that

$$\|f(p(s,t)) - \ell(s)\| < t \tag{5}$$

$$\operatorname{grad} f(p(s,t)) = \lambda(s,t) p(s,t) \tag{6}$$

$$\|p(s,t)\| > \frac{1}{t} \tag{7}$$

By Definable Choice [4] and by [9] ( see also Proposition 5.2 in [7] ) , after eventually shrinking the interval $J$ , we get $\varepsilon > 0$ such that

$$p(s,t) = a(s)t^\alpha + a_1(s,t)t^{\alpha_1}$$

$$f(p(s,t)) = \ell(s) + b(s)t^\beta + b_1(s,t)t^{\beta_1}$$

$$\lambda(s,t) = c(s)t^\gamma + c_1(s,t)t^{\gamma_1}$$

where $a, b, c \in C^1(J)$ and $a(s), b(s), c(s) \neq 0$ , $\forall s \in J$ , the functions $a_1, b_1, c_1 \in C^1(J \times (0,\varepsilon))$ are bounded on $J \times (0,\varepsilon)$ , and $\alpha < 0$ , $\alpha < \alpha_1$ , $\beta < \beta_1$ , $\gamma < \gamma_1$ . Moreover , by Lemma 2.5 , shrinking $J$ and reducing $\varepsilon$ if necessary , we may assume that $\frac{\partial a_1}{\partial s}, \frac{\partial b_1}{\partial s}, \frac{\partial c_1}{\partial s}$ are bounded on $J \times (0,\varepsilon)$ . By (5) , we have $\beta \geq 1$ . Therefore

$$\frac{\partial f(p(s,t))}{\partial s} \longrightarrow \ell'(s) \text{ when } t \longrightarrow 0 \quad .$$

Relation (6) implies

$$\frac{\partial f(p(s,t))}{\partial s} = \left\langle \frac{\partial p(s,t)}{\partial s}, \operatorname{grad} f(p(s,t)) \right\rangle = \overline{\lambda(s,t)} \left\langle \frac{\partial p(s,t)}{\partial s}, p(s,t) \right\rangle \quad .$$

Hence

$$\gamma + 2\alpha \leq 0 \quad . \tag{8}$$

On the other hand , by Lemma 2.6 , $t\dfrac{\partial f(p(s,t))}{\partial t} \longrightarrow 0$ when $t \longrightarrow 0$ . Relation (6) implies

$$t\frac{\partial f(p(s,t))}{\partial t} = t\overline{\lambda(s,t)} \left\langle \frac{\partial p(s,t)}{\partial t}, p(s,t) \right\rangle \longrightarrow 0 \text{ when } t \longrightarrow 0 \quad .$$

Hence , $1 + \gamma + (2\alpha - 1) > 0$ , i.e. $\gamma + 2\alpha > 0$ . This contradicts (8) . ◇

# 3 Examples and Remarks

**3.1** The Pfaffian functions ( see [6] for the definition ) , for example all functions $f \in \mathbf{R}[x_1, \ldots, x_n, \exp(x_1), \ldots, \exp(x_n)]$ , are shown to be definable in a suitable o–minimal structure . This is a consequence of a general result of Wilkie , see [13] .

**3.2 Remark .** (a) Let $g : \mathbf{R} \longrightarrow \mathbf{R}$ be defined by $g(x) := x^3$ . Then $g$ is a $C^0$ trivial fibration , but $\mathrm{B}_g = \Sigma_g = \{0\}$ .
(b) The bifurcation set of $f(x) = x \sin x$ is infinite . Obviously , $f$ is not definable in any o–minimal structure .
(c) In general , the inclusions $\Sigma_f \subseteq \mathrm{B}_f \subseteq \Sigma_f \cup S_f$ cannot be replaced by equalities , as the following examples show us .

**3.3 Example .** (a) Let $f : \mathbf{R} \longrightarrow \mathbf{R}$ be defined by $f(x) := x \exp x$ . Then $\Sigma_f = \{-1/e\}$ , $S_f = \{0\}$ and $\mathrm{B}_f = \{0, -1/e\}$ .
(b) Let $g : \mathbf{R}^2 \longrightarrow \mathbf{R}$ be defined by $g(x,y) := x^2 y^2 + 2xy$ . Then $\Sigma_g = \{0\}$ , $S_g = \{-1\}$ and $\mathrm{B}_g = \{0, -1\}$ .
(c) Let $h : \mathbf{R}^2 \longrightarrow \mathbf{R}$ be defined by $h(x,y) = y \exp(2x) + \exp x$ . Then $S_h = \{0\}$ , but $\mathrm{B}_h = \Sigma_h = \emptyset$ ( one can check that $H : \mathbf{R}^2 \longrightarrow h^{-1}(0) \times \mathbf{R}$ defined by $H(x,y) = (x, -\exp(-x), h(x,y))$ is a $C^\infty$ trivialization of $h$ ) .

Note that in examples (a) and (b) we have $\mathrm{B}_f = \Sigma_f \cup S_f$ , i.e. the equality can be attained .

**3.4 Remark .** Theorem 1.2 in [8] describe an approximation from above of the bifurcation sets of complex polynomial functions , using the Newton polyhedron at infinity. With the same proof as in [8] one can obtain a similar result for the case of *real* polynomial functions, and hence, the estimation from above becomes more effective .

**3.5 Remark .** If $f$ is a Pfaffian function , then using Proposition 1.6 and Khovanskii's theory on Fewnomials , see [6] , one can estimate from above the cardinalities of $S_f$ , $\Sigma_f$ and $\mathrm{B}_f$ . Note also that the conclusion of Proposition 1.6 is still true when $f : \mathbf{C}^n \longrightarrow \mathbf{C}$ is a complex polynomial , and hence , Khovanskii's theory can be applied to get a ( rough ) estimation of the cardinality of the bifurcation set in this case .

**Aknowledgements .** *We would like to thank The Fields Institute and the University of Toronto , where this note was written , for hospitality and support .*

Ta Lê Loi  
University of Dalat  
Department of Mathematics  
Dalat  
Vietnam  

Alexandru Zaharia  
Institute of Mathematics of  
the Romanian Academy  
*Mailing address ( until July* 1997 ) :  
The Fields Institute  
222 College Street  
Toronto  
Ontario, M5T 3J1  
Canada